\documentclass[12pt]{article}
\usepackage{latexsym,amsfonts,amsmath,amssymb}
\newtheorem{theorem}{Theorem}
\newtheorem{corollary}[theorem]{Corollary}
\newtheorem{sublemma}{Lemma}[theorem]
\newtheorem{lemma}[theorem]{Lemma}
\newtheorem{question}[theorem]{Question}
\newtheorem{observation}[theorem]{Observation}
\newtheorem{claim}[theorem]{Claim}
\newtheorem{subclaim}{Claim}[sublemma]
\newtheorem{conjecture}[theorem]{Conjecture}
\newtheorem{fact}[theorem]{Fact}
\newtheorem{definition}[theorem]{Definition}
\newtheorem{remark}[theorem]{Remark}
\newtheorem{example}[theorem]{Example}
\newtheorem{exercise}{Exercise}[section]
\def\Theorem #1.#2 #3\par{\setbox1=\hbox{#1}\ifdim\wd1=0pt
   \begin{theorem}{\rm #2} #3\end{theorem}\else
   \newtheorem{#1}[theorem]{#1}\begin{#1}\label{#1}{\rm #2} #3\end{#1}\fi}
\def\Corollary #1.#2 #3\par{\setbox1=\hbox{#1}\ifdim\wd1=0pt
   \begin{corollary}{\rm #2} #3\end{corollary}\else
   \newtheorem{#1}[theorem]{#1}\begin{#1}\label{#1}{\rm #2} #3\end{#1}\fi}
\def\Lemma #1.#2 #3\par{\setbox1=\hbox{#1}\ifdim\wd1=0pt
   \begin{lemma}{\rm #2} #3\end{lemma}\else
   \newtheorem{#1}[theorem]{#1}\begin{#1}\label{#1}{\rm #2} #3\end{#1}\fi}
\def\SubLemma #1.#2 #3\par{\setbox1=\hbox{#1}\ifdim\wd1=0pt
   \begin{sublemma}{\rm #2} #3\end{sublemma}\else
   \newtheorem{#1}{#1}[theorem]\begin{#1}\label{#1}{\rm #2} #3\end{#1}\fi}
\def\Question #1.#2 #3\par{\setbox1=\hbox{#1}\ifdim\wd1=0pt
   \begin{question}{\rm #2} #3\end{question}\else
   \newtheorem{#1}[theorem]{#1}\begin{#1}\label{#1}{\rm #2} #3\end{#1}\fi}
\def\Observation #1.#2 #3\par{\setbox1=\hbox{#1}\ifdim\wd1=0pt
   \begin{observation}{\rm #2} #3\end{observation}\else
   \newtheorem{#1}[theorem]{#1}\begin{#1}\label{#1}{\rm #2} #3\end{#1}\fi}
\def\Claim #1.#2 #3\par{\setbox1=\hbox{#1}\ifdim\wd1=0pt
   \begin{claim}{\rm #2} #3\end{claim}\else
   \newtheorem{#1}[theorem]{#1}\begin{#1}\label{#1}{\rm #2} #3\end{#1}\fi}
\def\SubClaim #1.#2 #3\par{\setbox1=\hbox{#1}\ifdim\wd1=0pt
   \begin{subclaim}{\rm #2} #3\end{subclaim}\else
   \newtheorem{#1}{#1}[sublemma]\begin{#1}\label{#1}{\rm #2} #3\end{#1}\fi}
\def\Conjecture #1.#2 #3\par{\setbox1=\hbox{#1}\ifdim\wd1=0pt
   \begin{conjecture}{\rm #2} #3\end{conjecture}\else
   \newtheorem{#1}[theorem]{#1}\begin{#1}\label{#1}{\rm #2} #3\end{#1}\fi}
\def\Fact #1.#2 #3\par{\setbox1=\hbox{#1}\ifdim\wd1=0pt
   \begin{fact}{\rm #2} #3\end{fact}\else
   \newtheorem{#1}[theorem]{#1}\begin{#1}\label{#1}{\rm #2} #3\end{#1}\fi}
\def\Definition #1.#2 #3\par{\setbox1=\hbox{#1}\ifdim\wd1=0pt
   \begin{definition}{\rm #2} {\rm #3}\end{definition}\else
   \newtheorem{#1}[theorem]{#1}\begin{#1}\label{#1}{\rm #2} {\rm #3}\end{#1}\fi}
\def\Remark #1.#2 #3\par{\setbox1=\hbox{#1}\ifdim\wd1=0pt
   \begin{remark}{\rm #2} {\rm #3}\end{remark}\else
   \newtheorem{#1}[theorem]{#1}\begin{#1}\label{#1}{\rm #2} {\rm #3}\end{#1}\fi}
\def\Example #1.#2 #3\par{\setbox1=\hbox{#1}\ifdim\wd1=0pt
   \begin{example}{\rm #2} #3\end{example}\else
   \newtheorem{#1}[theorem]{#1}\begin{#1}\label{#1}{\rm #2} #3\end{#1}\fi}
\def\Exercise #1.#2 #3\par{\setbox1=\hbox{#1}\ifdim\wd1=0pt
   {\footnotesize\begin{exercise}{\rm #2} {\rm #3}\end{exercise}}\else
   \newtheorem{#1}[section]{#1}{\footnotesize\begin{#1}\label{#1}{\rm #2} {\rm #3}\end{#1}}\fi}
\def\QuietTheorem #1.#2 #3\par{\setbox1=\hbox{#1}\ifdim\wd1=0pt\proclaim{Theorem {\rm #2}}{#3}\else\proclaim{#1 {\rm #2}}{#3}\fi}
\newcommand{\proclaim}[2]{\smallskip\noindent{\bf #1} {\sl#2}\par\smallskip}
\def\Proclaim #1.#2 #3\par{\proclaim{#1 {\rm #2}}{#3}}
\newenvironment{proof}{\noindent}{\kern2pt\QEDbox\par\bigskip}
\def\Proof#1: {\setbox1=\hbox{#1}\ifdim\wd1=0pt\begin{proof}{\bf Proof: }\else\medskip\begin{proof}{\bf #1: }\fi}
\newcommand{\QED}{\end{proof}}
\def\BF#1.{{\bf #1.}}
\def\Abstract #1\par{\begin{quotation}{\singlespaced\footnotesize{\noindent{\bf Abstract.~}#1}}\end{quotation}}
\def\Title #1\par{\title{#1}\maketitle}
\def\Author #1\par{\author{#1}}
\def\Acknowledgement#1\par{\thanks{#1}}
\def\Chapter #1\par{\chapter{#1}}
\def\Section #1\par{\section{#1}}
\def\QuietSection #1\par{\section*{#1}}
\def\SubSection #1\par{\subsection{#1}}
\def\SubSubSection #1\par{\subsubsection{#1}}
\def\MidTitle #1\par{\bigskip\goodbreak\centerline{\small\bf #1}\bigskip\noindent}

\newcommand{\singlespaced}{\baselineskip=15pt}
\def\bottomnote #1\par{{\renewcommand{\thefootnote}{}\footnotetext{#1}}}
\newcommand{\N}{{\mathbb N}}

\newfont{\msam}{msam10 at 12pt}

\newcommand{\jump}{{\!\triangledown}}
\newcommand{\Jump}{{\!\blacktriangledown}}
\def\ilt{<_{\infty}}
\def\ileq{\leq_{\infty}}
\def\iequiv{\equiv_{\infty}}
\newcommand{\Tequiv}{\equiv_T}

\newcommand{\Union}{\bigcup}

\newcommand{\trianglelt}{\lhd}

\newcommand{\smalllt}{\mathrel{\mathchoice{\raise2pt\hbox{$\scriptstyle<$}}{\raise1pt\hbox{$\scriptstyle<$}}{\scriptscriptstyle<}{\scriptscriptstyle<}}}
\newcommand{\smallleq}{\mathrel{\mathchoice{\raise2pt\hbox{$\scriptstyle\leq$}}{\raise1pt\hbox{$\scriptstyle\leq$}}{\scriptscriptstyle\leq}{\scriptscriptstyle\leq}}}

\newcommand{\UnderTilde}[1]{{\setbox1=\hbox{$#1$}\baselineskip=0pt\vtop{\hbox{$#1$}\hbox to\wd1{\hfil$\sim$\hfil}}}{}}
\newcommand{\Undertilde}[1]{{\setbox1=\hbox{$#1$}\baselineskip=0pt\vtop{\hbox{$#1$}\hbox to\wd1{\hfil$\scriptstyle\sim$\hfil}}}{}}
\newcommand{\undertilde}[1]{{\setbox1=\hbox{$#1$}\baselineskip=0pt\vtop{\hbox{$#1$}\hbox to\wd1{\hfil$\scriptscriptstyle\sim$\hfil}}}{}}
\newcommand{\UnderdTilde}[1]{{\setbox1=\hbox{$#1$}\baselineskip=0pt\vtop{\hbox{$#1$}\hbox to\wd1{\hfil$\approx$\hfil}}}{}}
\newcommand{\Underdtilde}[1]{{\setbox1=\hbox{$#1$}\baselineskip=0pt\vtop{\hbox{$#1$}\hbox to\wd1{\hfil\scriptsize$\approx$\hfil}}}{}}

\renewcommand{\th}{{\hbox{\scriptsize th}}}

\def\<#1>{\langle#1\rangle}

\newcommand{\QEDbox}{\fbox{}}

\newcommand{\WO}{\mathop{\hbox{\sc wo}}}

\newcommand{\factordiagramup}[6]{$$\begin{array}{ccc}
#1&\raise3pt\vbox{\hbox to60pt{\hfill$\scriptstyle
#2$\hfill}\vskip-6pt\hbox{$\vector(4,0){60}$}}&#3\\ \vbox
to30pt{}&\raise22pt\vtop{\hbox{$\vector(4,-3){60}$}\vskip-22pt\hbox
to60pt{\hfill$\scriptstyle #4\qquad$\hfill}}
     &\ \ \lower22pt\hbox{$\vector(0,3){45}$}\ {\scriptstyle #5}\\
\vbox to15pt{}&&#6\\
\end{array}$$}
\newcommand{\factordiagram}[6]{$$\begin{array}{ccc}
#1&&\\ \ \ \raise22pt\hbox{$\vector(0,-3){45}$}\ {\scriptstyle #2}
&\raise22pt\hbox{$\vector(2,-1){90}$}\raise5pt\llap{$\scriptstyle#3$\qquad\quad}&\vbox
to25pt{}\\ #4&\raise3pt\vbox{\hbox to90pt{\hfill$\scriptstyle
#5$\hfill}\vskip-6pt\hbox{$\vector(4,0){90}$}}&#6\\
\end{array}$$}
\newcommand{\cell}[1]{\boxit{\hbox to 17pt{\strut\hfil$#1$\hfil}}}
\newcommand{\head}[2]{\lower2pt\vbox{\hbox{\strut\footnotesize\it\hskip3pt#2}\boxit{\cell#1}}}
\newcommand{\boxit}[1]{\setbox4=\hbox{\kern2pt#1\kern2pt}\hbox{\vrule\vbox{\hrule\kern2pt\box4\kern2pt\hrule}\vrule}}
\newcommand{\Col}[3]{\hbox{\vbox{\baselineskip=0pt\parskip=0pt\cell#1\cell#2\cell#3}}}
\newcommand{\tapenames}{\raise 5pt\vbox to .7in{\hbox to .8in{\it\hfill input: \strut}\vfill\hbox to
.8in{\it\hfill scratch: \strut}\vfill\hbox to .8in{\it\hfill output: \strut}}}
\newcommand{\Head}[4]{\lower2pt\vbox{\hbox to25pt{\strut\footnotesize\it\hfill#4\hfill}\boxit{\Col#1#2#3}}}
\newcommand{\Dots}{\raise 5pt\vbox to .7in{\hbox{\ $\cdots$\strut}\vfill\hbox{\ $\cdots$\strut}\vfill\hbox{\
$\cdots$\strut}}}
\renewcommand{\dots}{\raise5pt\hbox{\ $\cdots$}}
\begin{document}
\author{Joel David Hamkins${}^1$\\
\normalsize\sc The City University of New York\\
{\footnotesize http://jdh.hamkins.org}}
\date{}

\Title Infinite time Turing machines

\Abstract Infinite time Turing machines extend the operation of ordinary Turing machines into transfinite
ordinal time.   By doing so, they provide a natural model of infinitary computability, a theoretical
setting for the analysis of the power and limitations of supertask algorithms.

Imagine that we have scaled the vertical asymptote of technological progress, as computers increase in
speed each day, and that we have in our possession an infinitely fast computer.  What would we do with it?
What {\it could} we do with it?  In order to make sense of these questions, we need a formal model of
infinite computability, a mathematical description of the computing machines and of the algorithms that
they will be able to carry out.  In this article, I introduce such a model of infinitary computability by
extending the classical Turing machine model into transfinite ordinal time; the resulting machines,
infinite time Turing machines, are able to carry out and complete algorithms involving infinitely many
computational steps. These computational steps proceed in time like the ordinal numbers, so that if the
computation does not halt at any of the finite stages 0, 1, 2, 3, $\ldots$, then it arrives at the first
infinite stage $\omega$, continuing with stages $\omega + 1$, $\omega + 2$, $\omega + 3$, $\ldots$ and so
on, eventually arriving at the second limit stage $\omega + \omega$, continuing with stages $\omega +
\omega + 1$, $\omega + \omega + 2$, $\ldots$ and so on through the ordinal numbers. The configuration of
the machine at each stage is determined by the earlier configurations and the operation of the fixed finite
program that the machine is running. These machines provide a robust model of infinitary computability, one
that has shed light on the power and limitations of supertask computation.

The essence of an infinitely fast computer, of course, is that it is able to complete infinitely many steps
of computation in a finite amount of time.  An infinitely fast computation, therefore, when measured by the
number of computational steps in it, can be viewed (perhaps paradoxically) as being infinitely long.  And
indeed, since the logical fundamental unit of a computation is the individual computation step, rather than
the length of time that such a step takes to carry out, it may be more natural to view infinitely fast
computations as being infinitely long.  Therefore, in this article let us focus on the infinitely long
computations, and take the term ``infinite time" computation, a term which might otherwise indicate a
painfully slow computation, to suggest instead a lightening fast computation that actually completes all
the steps in an infinite algorithm; the computational time is infinite when measured by the number of
computational steps.

\Section Supertasks

A {\it supertask} is a task involving infinitely many steps.  Because any kind of infinitary computation
will naturally involve algorithms having infinitely many computational steps---supertask algorithms---the
supertask concept seems to be central to any understanding of infinitary computability.  So before
introducing the infinite time Turing machine model in the next section, with its explicit supertask notions
of computability and decidability, let me first discuss a few supertask examples more generally.

Zeno of Elea (ca.  450 B.C.) was perhaps the first to grapple with the supertask concept in his famous
paradoxical argument that it is impossible to go from here to there.  Before arriving, one must first get
halfway there, and before that one must get halfway to the halfway point, and so on, {\it ad infinitum}.
Because one cannot accomplish these infinitely many tasks, Zeno argued, all motion is impossible.  Thus, he
takes the impossibility of completing a supertask as the foundation of his reductio.

In the twentieth century philosophical literature, the puzzles and paradoxes continue.  Take Thomson's lamp
(Thomson 1954), for example, which is on for $1/2$ minute, off for $1/4$ minute, on for $1/8$ minute, and
so on. After one minute (the sum of the series $1/2+1/4+1/8+\cdots$), is it on or off? The literature is
full of answers. Another toy example is the super-$\pi$ machine, which writes out the successive digits of
$\pi$ on a tape, the first in $1/2$ minute, the next in $1/4$ minute, and so on, so that all the digits are
written in a finite amount of time.  Because there can be no last or final step in such a process, Chihara
(1965) has criticized the completion of such an algorithm as unintelligible.

In a more entertaining example,\footnotemark[2] let's suppose that you have infinitely many one dollar
bills (numbered 1, 3, 5, $\ldots$) and upon entering a nefarious underground bar, you come upon the Devil
sitting at a table piled high with money.  You sit down, and the Devil explains to you that he has an
attachment to your particular bills and is willing to pay you a premium to buy them from you. Specifically,
he is willing to pay two dollars for each of your one-dollar bills.  To carry out the exchange, he proposes
an infinite series of transactions, in each of which he will hand over to you two dollars and take from you
one dollar. The first transaction will take $1/2$ hour, the second $1/4$ hour, the third $1/8$ hour, and so
on, so that after one hour the entire exchange will be complete.  The Devil takes a sip of whiskey while
you mull it over; should you accept his proposal? Perhaps you think you will become richer, or perhaps you
think with infinitely many bills it will make no difference? At the very least, you think, it will do no
harm, and so the contract is signed and the procedure begins.  How could the deal harm you?

It appears initially that you have made a good bargain, because at every step of the transaction, he pays
you two dollars and you give up only one.  The Devil is very particular, however, about the order in which
the bills are exchanged.  The contract stipulates that in each sub-transaction he buys from you your
lowest-numbered bill and pays you with higher-numbered bills.  Thus, on the first transaction he accepts
from you bill number 1, and pays you with bills numbered 2 and 4.  On the next transaction he buys from you
bill number 2 (which he had just paid you) and gives you bills numbered 6 and 8.  Next, he buys bill number
3 from you with bills 10 and 12, and so on.  When all the exchanges are completed, what do you discover?
You have no money left at all! The reason is that at the $n^\th$ exchange, the Devil took from you bill
number $n$, and never subsequently returned it to you.  So while it seemed as though you were becoming no
poorer with each exchange, in fact the final destination of every dollar bill in the transaction is under
the Devil's ownership.  The Devil is a shrewd banker, and you should have paid more attention to the
details of the supertask transaction to which you agreed.  And similarly, the point is that when we design
supertask algorithms to solve mathematical questions, we must take care not to make inadvertent assumptions
about what may be true only for finite algorithms.

Supertasks have also been studied by the physicists.  In classical mechanics under the Newtonian gravity
law (neglecting relativity), it is possible to arrange finitely many stars in orbiting pairs, each pair
orbiting the common center of mass of all the pairs, and a single tiny moon racing faster around, squeezing
just so between the dual stars so as to pick up speed with every such transaction.  Assuming point masses
(or collapsing stars to avoid collision), the arrangement can be made to lead by Newton's law of
gravitation to infinitely many revolutions in finite time.

Other supertasks reveal apparent violations of the conservation of energy in Newtonian physics: imagine
infinitely many billiard balls of successively diminishing size, converging to a point (see Laraudogoitia
1996).  The balls are initially at rest, but then the first is set rolling.  It collides with the next,
transferring all its energy, and that ball begins to roll.  Each ball in turn collides with the next,
transferring all its energy.  Because of the physical arrangements of the system, the motion disappears in
a sense into the singularity; after a finite amount of time all the collisions have taken place, and the
balls are at rest.  So we have described a physical process in which energy is conserved at each step, in
every interaction, but not overall through time.

Another arrangement has the balls spaced further and further out to infinity, pushing the singularity out
to the point at infinity.  The first ball knocks the second, which sends the next flying, and so on out to
infinity.  If the balls speed up sufficiently fast, under Newtonian rules it can be arranged that all the
motion is completed in a finite amount of time.  The interesting thing about this example and the previous
one is that time-symmetry allows them to run in reverse, with static configurations of balls suddenly
coming into motion without violating conservation of energy in any interaction.

More computationally significant supertasks have been proposed by physicists in the context of relativity
theory (Pitowsky 1990, Earman 1995, Earman, et. al. 1993, Hogarth 1992, Hogarth 1994).  Suppose that you
want to know the answer to some number theoretic conjecture, such as whether there are additional Fermat
primes (primes of the form $2^{2^n} + 1$), a conjecture that can be confirmed with a single numerical
example. The way to solve the problem is to board a rocket, while setting your graduate students to work on
earth looking for an example.  While you fly faster and faster around the earth, your graduate students,
and their graduate students and so on, continue the exhaustive search, with the agreement that if they ever
find an example, they will send a radio signal up to the rocket.  Meanwhile, by accelerating sufficiently
fast towards the speed of light (but never exceeding it), it is possible to arrange that the relativistic
time contraction on the rocket will mean that a finite amount of time on the rocket corresponds to an
infinite amount of time on the earth.   The point is therefore that by means of the communication between
the two reference frames, what corresponds to an infinite search can be completed in a finite amount of
time.  And it is precisely this ability to carry out such infinite searches---to go to the horizon and back
in a finite amount of time---that will make infinite time Turing machines so powerful.

With more complicated arrangements of rockets flying around rockets, one can solve more complicated number
theoretic questions.  Hogarth (1994) has explained how to determine in this way the truth of any arithmetic
statement.  More unusual relativistic Malament-Hogarth spacetimes can be (mathematically) constructed to
avoid the unpleasantness of unbounded acceleration required in these rocket examples (see also Hogarth
1992, Pitowsky 1990).

These computational examples speak to a widely held strengthening of Church's thesis, namely, the view that
the classical theory of computability has correctly captured the notion of what it means to be effectively
computable.\footnotemark[3]   Because the relativistic rocket examples provide algorithms for computing
functions, such as the halting problem, that are not computable by Turing machines, one can view them as
refuting this thesis. Supporters of this view emphasize that when thinking about what is in principle
computable, we must attend to the computational power available to us as a consequence of the fact that we
live in a relativistic or quantum-mechanical universe.  To ignore this power is to pretend that we live in
a Newtonian world.  Other (perhaps less successful) arguments against the thesis appeal to the random or
statistical nature of quantum mechanical phenomenon to mechanically produce non-computable bit-streams;
this is possible even in Newtonian physics:  a particle undergoing Brownian motion can be used to generate
a random bit stream that we have no reason to think is recursive.  Therefore, proponents argue, we have no
reason to believe the strong form of Church's thesis.

Apart from the question of what one can physically compute in this world, whether Newtonian or relativistic
or quantum-mechanical, a mathematician can be interested in what {\it in principle} a supertask can
accomplish. What are the mathematical limitations, as opposed to the physical inconveniences, of supertask
computation? Which mathematical questions can and which cannot be answered via supertasks?  To answer such
questions, one needs a precise mathematical model of the supertask machines and the algorithms that they
can carry out.  One can then closely analyze the model and attempt to place the class of decidable sets,
sets for which there is a computational procedure to decide whether a given object is a member of them or
not, within the known complexity hierarchies of descriptive set theory (which provide a way of measuring
the complexity of various sets of natural numbers or reals in terms of the complexity of their
definitions). The point is that such an abstract set-theoretic analysis can place severe bounds on the
class of sets that can be computable:  from below, we can hope to show that all sets at or below a given
level of complexity are computable; from above, we can hope to show that no sets above a certain level of
complexity are computable.  The proximity of these bounds to each other speaks to the quality of the
analysis, and precisely locates the supertask notion of computability within the larger (and intensely
studied) hierarchy of complexity.

Decades ago, Buchi (1962) and others initiated this kind of study at the lower levels of complexity, with
the study of $\omega$-automata and Buchi machines, involving automata and Turing machine computations of
length $\omega$ that accept or reject infinite input.  The editor of this volume, Jack Copeland, has been
particularly interested in the possibility of such machines completing their tasks in a finite amount of
time, by accelerating the computation in time along a convergent series (e.g. the first operation takes
$1/2$ second, the next $1/4$ second, and so on), and has put forth a philosophical analysis of such
accelerated Turing machines in various articles (see Copeland 1998a, 1998b, 2002).  Moving up in the
hierarchy, Gerald Sacks and many others founded the field of higher recursion theory, including
$\alpha$-recursion and $E$-recursion, a huge, deep body of work analyzing an abstract notion of computation
on infinite objects (see Sacks 1990).

In this article, I describe another model of infinitary computability, infinite time Turing machines, a
model which offers the strong computational power of higher recursion theory while remaining very close in
spirit to the computability concept of ordinary Turing machines.  Infinite time Turing machines were
originally defined by Jeff Kidder in 1989, and he and I worked out the early theory together while we were
graduate students at the University of California in Berkeley.  Six years later, Andy Lewis and I solved
some of the early questions and presented a complete introduction in a joint paper (appearing, finally, in
2000); in 1997 we solved Post's problem for supertasks (in press).  Benedikt Loewe (2001), Dan Seabold
(with myself, 2001) and especially Philip Welch (1999, 2000a, 2000b) have also made important
contributions.  The research is at an exciting early stage.

\Section Infinite time Turing machines

Infinite time Turing machines are obtained by extending the Turing machine concept to transfinite ordinal
time.  The idea is to allow a Turing machine to compute for infinitely many steps, preserving the
information produced while doing so, and allowing the machine to continue with more computation afterwards.
The individual computational steps will therefore be ordered like the ordinal numbers 0, 1, 2, 3, $\ldots$
$\omega$, $\omega+1$, $\omega+2$, $\ldots$ $\omega+\omega$, $\omega+\omega+1$, $\ldots$
$\omega+\omega+\omega$, $\ldots$ and so on. After each stage $\alpha$ of computation there is a unique next
stage $\alpha+1$, and these culminate in the limit ordinal stages, like $\omega$, $\omega+\omega$, and so
on.

So let me explain specifically how the machines work.  The machine hardware is identical to a classical
Turing machine, with a head moving back and forth reading and writing zeros and ones on a tape according to
the rigid instructions of a finite program, with finitely many states.  What is new is the transfinite
behavior of the machine, behavior providing a natural theory of computation on the reals that directly
generalizes the classical finite theory to the transfinite.

For convenience, the machines have three tapes-one for the input, one for scratch work and one for the
output-and the computation begins with the input written out on the input tape, with the head on the
left-most cell in the {\it start} state.

\begin{figure}[h]
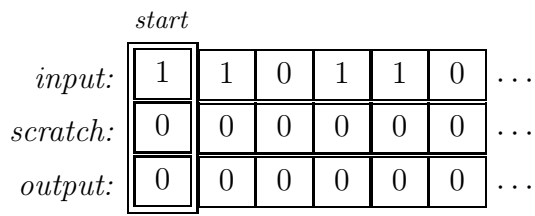

$$\tapenames\Head100{start}\Col100\Col000\Col100\Col100\Col000\Dots$$
\caption{An infinite time Turing machine---the computation begins}
\end{figure}

The successive steps of computation proceed in exactly the classical manner: the head reads the contents of
the cells on which it rests, reflects on its state and follows the instructions of the finite program it is
running.  Accordingly, it writes on the tapes, moves the head one cell to the left or the right or not at
all and switches to a new state.  Thus, the classical procedure determines the configuration of the
machine-the state of the machine, the position of the head and the contents of the tapes-at any stage
$\alpha + 1$, given the configuration at any stage $\alpha$.

We extend the computation into transfinite ordinal time by simply specifying the behavior of the machine at
the limit ordinal stages, such as $\omega$ or $\omega+\omega$.  When a classical Turing machine fails to
halt, it is often thought of in contemporary computability theory as a sort of failure; the result is
discarded even though the machine might have been writing some very interesting information on the tape,
such as all the theorems of mathematics, for example, or the members of some other computably enumerable
set.  With infinite time Turing machines, however, we hope to preserve this information by taking some kind
of limit of the earlier configurations and continuing the computation transfinitely.

\begin{figure}[h]
$$\tapenames\Head101{limit}\Col111\Col010\Col101\Col001\Col011\Dots$$
\caption{The limit configuration}
\end{figure}

Specifically, at any limit ordinal stage, the head resets to the left-most cell; the machine is placed in
the special {\it limit} state, just another of the finitely many states; and the values in the cells of the
tapes are updated by computing a kind of limit of the previous values that cell has displayed.
Specifically, if the values in a cell have stabilized before a limit stage, then the limit value displayed
by that cell at the limit stage will be this stabilized value; and otherwise, when the cell's value has
alternated from 0 to 1 and back again unboundedly often before a limit stage, then the limit value is set
to 1.  This limit value is equivalent, by definition, to computing for each cell what is called in
mathematics the {\it limit-supremum} (the ``$\limsup$") of the previous values displayed in that cell. (The
model is robust in that if we use a limit value of 0 in such cases, corresponding instead to the {\it
limit-infimum}, an equivalent model of computability is obtained.) With the limit stage configuration thus
completely specified, the machine simply continues computing.  If after some amount of time the {\it halt}
state is reached, the machine gives as output whatever is written on the output tape.

\begin{figure}[h]
$$\tapenames\Col101\Col011\Col000\Head001{halt}\Col110\Col011\Dots$$
\caption{The computation ends, giving output $110101\cdots$}
\end{figure}

Because there seems to be no need to limit ourselves to finite input and output---the machines have plenty
of time to consult the entire input tape and to write on the entire output tape before halting---the
natural context for these machines is Cantor Space $2^\omega$, the space of infinite binary sequences.  For
our purposes here, let's refer to such infinite binary sequences as real numbers.\footnotemark[4]

In this way, the machines provide a model of computation on the real numbers.  Specifically, if $p$ is a
infinite time Turing machine program, then the function that $p$ computes is determined by simply running
an infinite time Turing machine with program $p$ on whatever input is given, and giving as output the
output of that computation.  Such functions are said to be {\it infinite time computable}, or {\it
supertask computable}. (Please note that since some computations never halt, even in the infinite time
context, and so the domain of such a function need not include all real numbers.)  A set of reals $A$ is
{\it infinite time decidable}, or {\it supertask decidable}, when its characteristic function, the function
with value 1 for inputs in $A$ and 0 for inputs not in $A$, is supertask computable.  The set $A$ is
infinite time {\it semi-decidable} when there is a program that halts with output 1 on exactly the inputs
that are in $A$. (Thus, the semi-decidable sets correspond in the classical theory to the computably
enumerable sets, though since here we have sets of reals, one hesitates to describe them as enumerable.)
Since it is an easy matter to change any output value to 1, the semi-decidable sets are exactly the sets
that are the domain of a computable function, just as in the classical theory.

The first mathematical observation to make is that one need only consider countable ordinals (these are the
ordinals that can be put into a one-to-one correspondence with a set of natural numbers, as opposed to the
uncountable ordinals, which cannot) when analyzing the computations of infinite time Turing machines.  The
fact is that every supertask computation either halts or repeats in countably many steps.  Let me explain
in some mathematical detail why this is so.  Suppose that a supertask computation has not halted in
countably many steps; I will argue that it is caught in an endlessly repeating loop.  By assumption, the
computation proceeded through all the countable ordinals without halting, and finds itself at stage
$\omega_1$, the first uncountable ordinal (the smallest ordinal having cardinality $\aleph_1$), without
having yet halted. Consider the complete configuration of the machine---the position of the head, the state
and the contents of the tape---at this stage.  I claim that this configuration must have occurred earlier,
and the machine is caught in an endlessly repeating loop.  Since $\omega_1$ is a limit ordinal, we know
that the head is on the left-most cell in the {\it limit} state, as it is at any limit ordinal stage.
Further, any cell with a value of 0 at stage $\omega_1$ must have been 0 from some countable ordinal on,
and since there are only countably many cells, there must be a countable ordinal stage $\alpha_0$ by which
all of those 0 valued cells have stabilized. Let $\alpha_1$ be a larger countable ordinal by which all the
other cells have shown a value of 1 at least once, and $\alpha_2$ be a still larger value by which all the
remaining cells have shown a value of 1 once again, and so on. Let $\alpha_\omega$ denote the first stage
after all the stages an (this is still a countable ordinal). I claim that the configuration at stage
$\alpha_\omega$ is the same as it is at stage $\omega_1$: since it is a limit ordinal, the head is on the
first cell and in the {\it limit} state; and by construction the contents of each cell are computed in the
same way that they are at stage $\omega_1$. Further, since beyond $\alpha_0$ the only cells that change are
the ones that will change unboundedly often, by design, it follows that limits of this configuration are
the very same configuration again, and the machine is caught in an endlessly repeating loop.  (Please note
that in the supertask context, a computation that merely repeats a complete machine configuration is not
necessarily caught in an endlessly repeating loop, for after $\omega$ many repetitions, the limit
configuration may allow it to escape. This is why one must consider the limits of the repeating
configuration in the preceding argument.  An example of this phenomenon is provided by the machine which
does nothing at all except halt when it is in the {\it limit} state; this machine repeats its initial
configuration many times, yet still escapes this repeating loop and halts at $\omega$.)  The conclusion is
that any computation that fails to halt in countably many steps is doomed to repeat itself endlessly as the
ordinals fall through the transfinite hourglass.

\Section How powerful are the machines?

One naturally wants to understand the power of the new machines.  Since, of course, one needn't take
advantage of the new powers the machine offers, any function that is computable by a classical Turing
machine is computable by an infinite time Turing machine.  But further, I claim that the infinite time
Turing machines are strictly more powerful than their classical counterparts.  Specifically, the halting
problem for classical Turing machines---the question of whether a given program $p$ halts on a given input
n in finitely many steps, a problem well-known to be undecidable by classical Turing machines---is
decidable in $\omega$ many steps by infinite time Turing machines.  To see this, one programs an infinite
time Turing machine to simulate the operation of $p$ on $n$, and if the simulated computation ever halts
our algorithm gives the output that yes, indeed, the computation did halt.  Otherwise, the {\it limit}
state will be attained, and when this occurs the machine will know that the simulated computation failed to
halt; so it outputs the answer that no, the computation did not halt.  In this algorithm, it is as though
the {\it limit} state exists specifically to tell the machine that the simulated computation did not yet
halt.

This argument generalizes to show that the infinite time Turing machines can decide membership in any given
computably enumerable set $A$ in $\omega$ many steps.  On input $n$, one simulates the classical
enumeration of $A$; if $n$ appears in the enumeration, give the output that yes, $n$ is in $A$.  Otherwise,
when the {\it limit} state appears, give the output that no, $n$ is not in $A$.

The power of infinite time Turing machines, though, far exceeds the classical halting problem and the
computably enumerable sets.  In fact, any question of first order number theory is supertask decidable.
With an infinite time Turing machine, one could solve the prime pairs conjecture (which asserts that there
are infinitely many prime pairs---pairs of primes, such as $\<11, 13>$ or $\<107, 109>$, which differ by
two), for example, and the question of whether there are infinitely many Fermat primes (primes of the form
$2^{2^n}+1$) and so on. With an infinite time Turing machine, one can decide the truth of any statement
that is expressible in first order number theory (this includes essentially all of elementary number
theory, all of finite combinatorics, all of finite graph theory, and so on).  To decide a question of the
form $\exists n\, \varphi(n,x)$, for example, where $n$ ranges over the natural numbers, one can simply try
out all the possible values of $n$ in turn.  One either finds a witness $n$ or else knows at the limit that
there is no such witness, and in this way decides whether $\exists n\, \varphi(n,x)$ is true.

By analyzing this idea precisely, once concludes by induction on the complexity of the statement that any
first order number theoretic question is decidable with only a finite number of limits, that is, before
stage $\omega^2$ (the limit of the sequence $\omega$, $\omega\cdot 2$, $\omega\cdot 3$, and so on). In
fact, it is possible to prove that the class of sets that are decidable in time uniformly before $\omega^2$
(meaning that for some $n$ the program halts on all input before stage $\omega\cdot n$) is exactly the
class of arithmetic sets, the sets of reals that are definable by a statement using first order quantifiers
over the natural numbers (see Hamkins, et. al. 2000, Theorem 2.6).

Infinite time Turing machines, however, have a power that transcends even the full first order theory.  In
fact, their power reaches into the projective hierarchy, in the second order theory (the projective sets
are those that are definable using quantifiers not only over the natural numbers, but also over sets of
natural numbers, or over the real numbers).  One particular set, known as $\WO$ (for well-order, since it
consists of the real numbers that code a well-ordered relation on the natural numbers, a relation having no
infinite descending sequence), serves as a familiar sentinel or gate-keeper to these higher levels of
complexity.  Specifically, the set $\WO$ exists at what is called the $\Pi^1_1$ level of complexity, a
level that subsumes all of the arithmetic hierarchy of first order number theory, as well as the
hyperarithmetic sets (a kind a transfinite generalization of the arithmetic sets), and it is complete for
this level of complexity in the sense that every other $\Pi^1_1$ set reduces to it.  In this sense, $\WO$
is the hardest possible $\Pi^1_1$ set, and so to show that $\WO$ is supertask decidable, as I will do in
the next paragraph, is to show that the infinite time Turing machines have power completely transcending
the classical computability theory.

With the reader's permission, let me sketch the argument showing that $\WO$ is decidable, using what Andy
Lewis and I (2000) have called the ``count-through" argument.  The set $\WO$ consists of all reals coding a
well-order relation on a subset of the natural numbers.  (A real $x$ codes a relation $\trianglelt$ on $\N$
in the sense that $i \trianglelt j$ if and only if the $\<i,j>^\th$ bit of $x$ is 1, using the G\"odel
pairing function $\<i,j>$, which associates pairs of natural numbers with individual natural numbers in a
one-to-one correspondence. This function allows a binary relation to be coded with a binary sequence.) We
would like to describe a supertask algorithm which on input $x$ decides whether $x$ codes a well-ordered
relation $\trianglelt$ on a set of natural numbers or not, that is, whether $x \in \WO$ or not.  In
$\omega$ many steps, it is easy to check whether $x$ codes a linear order: this amounts merely to checking
that the relation $\trianglelt$ coded by $x$ is transitive, irreflexive and connected.  That is, the
machine must check that whenever $i \trianglelt j$ and $j \trianglelt k$ then also $i \trianglelt k$, that
$i \trianglelt i$ never holds, and that for distinct $i$ and $j$ either $i \trianglelt j$ or $j \trianglelt
i$. All these requirements can be enumerated and checked in $\omega$ many steps. Next, the algorithm will
attempt to find the least element in the field of the relation $\trianglelt$. This can be done by keeping a
current-best-guess on the scratch tape and systematically looking for better guesses, whenever a new
smaller element is found. When such a better guess is found, it replaces the current guess on the scratch
tape, and a special flag cell is flashed on and then off again. At the limit, if the flag cell is on, it
means that infinitely often the guess was changed, and so the relation has an infinite descending sequence.
Thus, in this case the input is definitely not a well order and the computation can halt with a negative
output.  Conversely, if the flag cell is off, it means that the guess was only changed finitely often, and
the machine has successfully found the $\trianglelt$-least element. The algorithm now proceeds to erase all
mention of this element from the field of the relation $\trianglelt$, by turning the bits $\<i,j>$ of $x$
to 0 that mention it. This produces a new smaller relation, and the algorithm proceeds to find the least
element of it.  In this way, the relation $\trianglelt$ is eventually erased from the bottom as the
computation proceeds.  If the relation is not a well order, eventually the algorithm will erase the largest
initial segment of it that is a well order, and then discover that there is no least element remaining, and
reject the input.  If the relation is a well order, then the algorithm will eventually erase the entire
field, and recognize that it has done so, and accept the input as a well order. So in any case the
algorithm can recognize whether $x$ codes a well order or not, and so $\WO$ is supertask decidable.  It
follows, as I explained, that every set at or below the $\Pi^1_1$ level of complexity is also decidable.
Because of this, the power of infinite time Turing machines is seen to reach into second order number
theory, completely transcending the classical theory.\footnotemark[5]

Much of the classical computability theory generalizes to the supertask context of infinite time Turing
machines.  For example, the classical $s$-$m$-$n$ theorem and the Recursion Theorem, for those who are
familiar with them, go through with virtually identical proofs.  But some other classical results, even
very elementary ones, do not generalize.  One surprising result, for example, is the fact that there is a
non-computable function whose graph is semi-decidable.  This follows from what Any Lewis and I (in our 2000
paper) have called the ``Lost Melody Theorem," which asserts the existence of a real $c$ which is not
writable, so that it is not the output of any computation on input 0, but such that the set $\{ c \}$ is
decidable. Imagine the real $c$ as the melody that you can recognize when someone sings it, but you cannot
sing it on your own; so you can say yes or no whether a given object is c or not, thereby deciding the set
$\{ c \}$, but you cannot create $c$ from scratch.  If $c$ is such a lost melody real, then let $f$ be the
function defined by $f(x) = c$.  This function is total (defined on all input) and constant (has the same
value for all input). The graph of the function consists of all pairs $(x,f(x))$, that is, of all pairs
$(x,y)$ for which $y$ is $c$, and this is decidable.  But the function is not computable because $c$ is not
writable.  To summarize, the example shows that in the supertask context we have situation completely
contrary to the classical computability theory: a non-computable function whose graph is semi-decidable;
indeed, the example is all the more surprising because graph of the function $f$ here is actually
decidable, not merely semi-decidable, and what is more, the function is total and constant!  There is
nothing like this in the classical computability theory.

\Section How long are the computations?

One naturally wants to understand how long a supertask computation can be.  Therefore, I define an ordinal
$\alpha$ to be {\it clockable} if there is a computation on input 0 that takes exactly $\alpha$ many steps
to complete (so that the $\alpha^\th$ step of computation is the act of moving to the {\it halt} state).
Such a computation is a clock of sorts, a way to count exactly up to $\alpha$.  One sets the machine
running, and when it halts it is as though an alarm clock has gone off:  exactly $\alpha$ many steps have
been computed. Such clocks are especially useful when it is desired to run some other computation for a
given number of steps; one simulates the clock and the other computation simultaneously, and when the clock
stops, the algorithm has computed the desired number of steps.

It is very easy to see that any finite number $n$ is clockable; one can simply have a machine cycle through
$n$ states and then halt.  The ordinal $\omega$ is clockable, by the machine that halts when it encounters
the {\it limit} state.  Similar ideas show that if $\alpha$ is clockable, then so is $\alpha + n$ and
$\alpha + \omega$. From this, it follows that every ordinal up to $\omega^2$ is clockable.  The ordinal
$\omega^2$ itself is clockable: one can recognize it as the first limit of limit ordinals, by flashing a
cell on and then off again every time the {\it limit} state is encountered; the ordinal $\omega^2$ will be
first time this cell is on at a limit stage. Going beyond this, it is easy to see that if $\alpha$ and
$\beta$ are clockable, then so are $\alpha + \beta$ and $\alpha\beta$. Students might enjoy finding
algorithms to clock specific ordinals, such as $\omega^{\omega^2} + \omega^4$, and I can recommend this as
a way to help them understand the ordinals more deeply.

Many readers will have guessed that the analysis extends much further, to much larger ordinals.  In fact,
any recursive ordinal is clockable.  (A {\it recursive ordinal} is one whose order relation is isomorphic
to recursive relation on the natural numbers.  The supremum of the recursive ordinals, of which there are
only countably many, is denoted $\omega_1^{CK}$, named for Church and Kleene.)  This can be seen by
optimizing the count-through argument I gave earlier when proving that $\WO$ is supertask decidable.
Specifically, after writing on the tape in $\omega$ many steps a real that codes a given recursive ordinal,
one proceeds to count through this ordinal by guessing the $\omega$ many least elements of the relation
(while simultaneously erasing the previous guesses), rather than handling only one guess with each limit.
In this way, each block of $\omega$ many steps of the algorithm will erase $\omega$ many elements from the
field of the relation, and the length of time the count-through procedure takes is more closely connected
with the ordinal being counted.

Some have been surprised that the clockable ordinals extend beyond the recursive ordinals, but in fact they
extend well beyond the recursive ordinals.  To see at least the beginnings of this, let me show that the
ordinal $\omega_1^{CK}+\omega$ is clockable.  By Feferman and Spector (1962), there is a recursive relation
whose well-founded part has order type $\omega_1^{CK}$.  Consider the supertask algorithm that writes this
relation on the tape and then attempts to count through it.  By stage $\omega_1^{CK}$ the non-well-founded
part of the relation will have been reached, but it takes the algorithm an additional $\omega$ many steps
to realize this.  So it can halt at stage $\omega_1^{CK}  + \omega$.

One is left to wonder, is $\omega_1^{CK}$ itself clockable? More generally, Are there gaps in the clockable
ordinals? After all, if a child can count to twenty-seven, then one might expect the child also to be able
to count to any smaller number, such as nineteen.\footnotemark[6]   The question is whether we expect the
same to be true for infinite time Turing machines.

The surprising fact is that there are gaps in the clockable ordinals.  Let me give a brief sketch of this
argument:  Consider the algorithm which simultaneously simulates all programs on input 0, recording which
have halted; when a stage is found at which no programs halt, then halt.  By construction, this algorithm
verifies that a supertask algorithm can halt above a non-clockable ordinal, and so there are gaps in the
clockable ordinals.  The argument can be modified to show that the next gap above any clockable ordinal has
length $\omega$.  Other arguments establish that complicated behavior can occur at limits of gaps, because
the lengths of the gaps are unbounded in the clockable ordinals.

An interesting question for future mathematical research would be to understand the structure of the
clockable ordinals.  For example, it is now known that no admissible ordinal can be clockable, and so the
first gap occurs exactly at $\omega_1^{CK}$.  Which is the next non-clockable ordinal?  Can we find an
independent characterization of the clockable ordinals? They do seem to arise very naturally in this
theory.

There is another way for infinite time Turing machines to operate as clocks, and this is by carefully using
the count-through argument that I mentioned earlier.  To assist with this analysis, we define that a real
is {\it writable} if it is the output of a supertask computation on input 0.  In analogy with the recursive
ordinals, we define that an ordinal is writable if it is coded by a writable real.  It is easy to see that
there are no gaps in the writable ordinals, because if one can write down real coding $\alpha$, it is an
easy matter to write down from this a real coding any particular $\beta < \alpha$.  Andy Lewis and I (2000)
proved that the order-types of the set of clockable ordinals and the set of writable ordinals are the same,
but the question was left open as to whether these two classes of ordinals had the same supremum, that is,
whether the two sets converged upwards to the same ordinal.  This was solved affirmatively by Philip Welch
(2000a), allowing Andy Lewis and I to greatly simplify arguments in our paper on Post's Problem.

\Section The supertask halting problems and jump operators

Before discussing the supertask notions (and at the suggestion of the editor of this volume), let me spend
several paragraphs giving an introduction to the elements of classical degree theory.  Once having
specified the operation of the classical Turing machines, one knows what it means for a function on the
natural numbers to be computable and for a set of natural numbers to be decidable or computably enumerable
(also known as recursively enumerable).  The halting problem, the set of programs halting on a given input,
then appears as a canonical example of a set that is computably enumerable but not decidable; indeed, it
may be first concrete set many encounter that they know definitively to be undecidable.  The natural
temptation arises to ask:  Suppose you had some way of answering the halting problem, what then? With this
set in hand, what could you compute from it?  To answer, one develops the theory of oracle computation.
Specifically, for a fixed set of natural numbers $A$, known as the oracle, one equips the Turing machine
with an extra tape, the oracle tape, on which the membership information of $A$ is explicitly written (the
$n^\th$ cell of the oracle tape is 1 when $n$ is in $A$, otherwise 0).  The functions that are computable
with this expanded device are the functions that are computable relative to $A$, and one has a notion of a
set being decidable or computably enumerable relative to $A$.

The notion of relative computability reveals very quickly a vast mathematical structure lurking behind the
scene: the Turing degrees.  Specifically, one defines that $A\leq_T  B$ (pronounced ``$A$ is Turing
reducible to $B$") when $A$ is decidable relative to oracle $B$, and this relation provides a way of
measuring the relative complexity of sets; if $A \leq_T  B$, then $B$ is at least as complicated as $A$.
The two sets are Turing equivalent, written $A \Tequiv  B$, when each oracle is computable relative to the
other, that is, when $A \leq_T B$ and $B \leq_T A$.  If $A$ is any set, then the degree of $A$ is the
equivalence class of $A$ in this relation, the collection of all $B$ such that $A \Tequiv  B$.  Any such
$B$ has the same information content as $A$, because with $B$ we may reconstruct $A$ and vice versa.  The
strict relation $A <_T B$ is defined to hold when $A \leq_T B$ but not $A \Tequiv  B$.  Two sets are
incomparable (written $A \perp_T B$) when neither is computable from the other.

The undecidability of the halting problem $H$ means that $0 <_T H$, where $0$ refers to the simplest
possible set, the empty set, whose degree is the collection of all decidable sets.  It is relatively easy
to see that every computably enumerable set is decidable from $H$, and so the question remains, known as
Post's Problem, whether there are any computably enumerable sets strictly between $0$ and $H$?  That is,
given the examples of the decidable sets and the various forms of the halting problem, which are all Turing
equivalent, one might think that the only way to show that a given problem is undecidable is to reduce the
halting problem to it. This would mean that every computably enumerable set is either decidable or
equivalent to the halting problem, and there is nothing strictly between them.  Post's Problem was solved
by Friedburg and Muchnik in 1956, however, with the opposite answer: with their finite injury priority
argument method, one can construct a rich diversity of sets intermediate between $0$ and $H$ (see Soare
1987). We now know that there is a rich structure of computably enumerable degrees; every countable partial
order, for example, finds an isomorphic copy in these degrees.

It is surprisingly easy to lift the picture that we have of the degrees between $0$ and $H$ to higher
levels of the hierarchy, for the simple reason that most computability arguments relativize to the presence
of a fixed oracle $A$.   Thus, for example, the relativized halting problem $H^A$, the set of all programs
that halt on a given input with oracle $A$, is computably enumerable but not decidable relative to $A$.
This observation provides a method of jumping higher in the degree structure:  the Turing {\it jump} of a
set $A$, written $A'$ and pronounced ``$A$-jump," is the halting problem relative to $A$, and by
relativizing the undecidability of the halting problem, one knows that $A <_T A'$.   In this terminology,
the classical halting problem $H$ becomes $0'$, pronounced ``$0$-jump".  Faced with any operator, of
course, one is tempted to iterate it, and doing so here produces an increasing sequence of degrees reaching
upwards into the Turing hierarchy:  $0 <_T 0' <_T 0''  <_T 0'''  <_T \cdots <_T 0^{(n)} <_T \cdots$ and so
on. One can ask all kinds of questions about the structure of the Turing degrees:  Is every degree above
$0'$ the jump of another degree? Yes.  If the jumps of two sets are Turing equivalent, are the sets
themselves equivalent? No.  Are there undecidable degrees below $0'$ whose jump is still only $0'$? Yes.
Are there minimal non-computable degrees? Yes.  Are the computably enumerable degrees dense, in the sense
that whenever one is strictly less than another, then there is a third strictly in between?  Yes. (For a
textbook account of all these answers and much more, consult Soare 1987).

Perhaps the broad philosophical view here is that each Turing degree represents a possible amount of
information---this being precisely what is preserved by the Turing equivalence relation---and if we are
interested in the structure of information, which I assume most readers of this volume are, it behooves us
to understand the structure of the Turing degrees.

This concludes my introduction to the theory of the classical Turing degrees; let me now turn to the
corresponding supertask model.  Every new model of computation naturally also provides a corresponding new
halting problem, the question whether a given computation halts.  Thus, in the supertask context of
infinite time Turing machines, we have the supertask halting problem.  The full supertask halting problem
$H$ is the question of whether a given infinite time Turing machine program $p$ halts on a given real input
$x$. A special case of this problem is the weak halting problem $h$, which is the question of whether a
given program $p$ halts on the trivial input $0$.  The analogues of these two problems in the classical
computability context are Turing equivalent and there is no big fuss made there distinguishing them, but
here in the supertask context the situation is somewhat subtler.

The classical arguments directly generalize to show that the halting problems $h$ and $H$ are
semi-decidable but not decidable.  First, for semi-decidability, the point is that given a program $p$ and
input $x$ (or input $0$), one can simply simulate $p$ on $x$ to see if it halts.  If it does, output the
answer that yes, it halted; otherwise, keep simulating.  Second, for undecidability, let me give in the
case of $H$ the classical diagonalization argument.  Suppose that $H$ were decidable by an infinite time
Turing machine, and consider the following algorithm: on input $p$, use the decision procedure for $H$ to
determine whether $p$ halts on input $p$; if it does, then go into an infinite loop and never halt; if it
does not, then halt immediately. This algorithm is infinite time computable by some program $q$.  But
observe that, by design, we have achieved the ubiquitous diagonal contradiction, namely, the program $q$
halts on input $q$ if and only if $q$ does not halt on input $q$.  So $H$ cannot be decidable.  The
classical arguments also generalize to show that the weak halting problem $h$ is also undecidable.

Once one has introduced the halting problems, it is again natural to wonder what one could compute if one
had some means to solve these problems, and one is led to the notion of oracle computation.  There are,
however, two natural types of oracles to use in the infinite time Turing machine context.  On the one hand,
one can use an individual real number as an oracle just as one does in the classical context, by simply
adding an oracle tape containing this real number, and allowing the machine to access this tape during the
computation.  This corresponds exactly to adding an extra input tape and thinking of the oracle as a fixed
additional input.

But ultimately one realizes this is not the right type of oracle to consider, since an oracle is more
properly the same type of object as one that might be decidable or semi-decidable, namely, a set of real
numbers, not a mere individual real number.  That is, with supertasks we have made the move to a
higher-type theory of computability-with functions defined on sets of real numbers, rather than on sets of
natural numbers-and we naturally want to develop our theory of relative computability at this level as
well.  Thus, we want somehow to allow a set of real numbers, such as the halting problem $H$, to become an
oracle for the machines.  And since such a set could be uncountable, and the set $H$ in particular is
definitely uncountable, we can't expect to be able to write out the entire contents of the oracle on an
extra tape, as in the classical theory.  Instead, we provide an oracle model of relative computability by
which the machine can make membership queries of the oracle.  Specifically, for a fixed oracle set of reals
$A$, we equip an infinite time Turing machine with an initially blank oracle tape on which the machine can
read and write.  The machine makes a membership query of $A$ by writing out a real number $a$ on the oracle
tape, and then switching to a special query state; the result, however, is that the machine moves actually
to the yes or no state, depending on whether $a\in A$ or not, and the algorithm has completed its query of
$A$. In this way, the machine is able to ask, of any real $a$ that it is capable of producing, whether
$a\in A$ or not, and use the answer as the basis for more computation.  This query model of oracle
computation has proven robust. In the lower-type case of oracle sets of natural numbers, it is
computationally equivalent to the usual oracle tape model in both the classical theory and the supertask
theory here.  And in the higher-type case of sets of reals, it closely follows the well-known definition in
set theory of G\"odel's constructible universe $L[A]$, relative to the predicate $A$, since when building
this model, one is allowed to apply the predicate only to previously constructed objects, just as our
oracle machines are only allowed to make queries of the form ``is $a\in A$?" for real numbers $a$ that they
are able to produce for the purpose. Queries, for example, of the form ``is there a real $y$ in the oracle
having such-and-such property?" do not have the required form.

From the notion of oracle computation, of course, just as in the classical theory, one obtains a notion of
relative computability.   Specifically, one set $A$ is computable {\it relative} to another set $B$, if
having $B$ as an oracle allows one to decide membership in $A$, that is, if the characteristic of $A$ is
computable with oracle $B$.  When this occurs, we write $A \ileq  B$, indicating that $B$ is at least as
complicated as $A$, and we say that $A$ reduces to $B$.  From this notion of relative computability, we can
then define when two oracles are {\it equivalent}:  $A$ is equivalent to $B$, written $A \iequiv  B$, if
and only if each is computable from the other, that is, $A \ileq B$ and $B \ileq A$; the strict relation $A
\ilt B$ is defined to hold if and only if $A \ileq B$ but not $A \iequiv  B$. The equivalence classes of
the equivalence relation $\iequiv$ form the infinite time Turing {\it degrees}, namely, two sets have the
same degree when they are equivalent to each other. When two sets are equivalent, they have the same
information content, as far as supertask computation is concerned, because from either of them one can
consruct the other.  It is therefore in the structure of these degrees that nature of supertask computation
is to be revealed; it is here that we will discover both the power and the limitations of supertask
computation.

The two halting problems give rise, of course, to two jump operators.  Specifically, for any set $A$ we
have the halting problem $H$ relative to $A$, the question of whether a given program $p$ halts on a given
input $x$ using oracle $A$.  This problem is denoted by $A^\Jump$, pronounced $A$-jump, in analogy with the
Turing jump. We also have the weak jump $A^\jump$, obtained by adding to $A$ the weak halting problem $h$
relative to oracle $A$. The Jump Theorem asserts that for any set, $A \ilt A^\jump \ilt A^\Jump$, and it is
because of this theorem that we know that the two halting problems $h$ and $H$ are not equivalent.  This is
because $h$ is equivalent to $0^\jump$ and $H$ is equivalent to $0^\Jump$, and so since $0^\jump \ilt
0^\Jump$, it follows that $h \ilt H$. That is, $h$ is computable from $H$ but not vice versa.  The Jump
Theorem is generalized by the absorption property of the jump operators, namely, the fact that the strong
jump absorbs applications of the weak jump, expressed by the equivalence $A^{\jump\Jump} \iequiv A^\Jump$.
This theorem shows that in a sense, the strong jump $\Jump$ jumps much higher than the weak jump $\jump$,
as it means that $A^\Jump$ is much higher in the degree structure than $A^\jump$. In particular, it follows
from the absorption property that $A^\jump$ cannot be equivalent to $A^\Jump$, since in that case
$A^{\jump\Jump}$ would be equivalent to $A^{\Jump\Jump}$, which by the Jump Theorem is strictly more
complicated than $A^\Jump$. All these results are discussed at length in my paper with Andy Lewis (2000).

\Section Post's Problem for Supertasks

Post's problem, the question in classical computability theory of whether there are non-decidable
semi-decidable degrees other than the halting problem, has a direct supertask analogue: Are there
intermediate semi-decidable degrees between $0$ and the jump $0^\jump$?

The answer is delicately mixed.  On the one hand, contrasting sharply with the classical theory in the
context of degrees in the real numbers, we have a negative answer.  The fact is that there are no real
numbers $z$ strictly in between 0 and the jump $0^\jump$.  Thus, while over the past fifty years the
classical theory has supported an active, vibrant study of the computably enumerable degrees and indeed, of
all the degrees between 0 and the Turing jump $0'$, my theorem with Andy Lewis shows that in the context of
supertask degrees of real numbers between 0 and the jump $0^\jump$, there is simply nothing to study.

Since this result is one of the strongest departures from the classical theory that has yet arisen, let me
briefly sketch the proof.  Suppose that $z$ is a real number and $0 \ileq z \ileq 0^\jump$. So $z$ is the
output of program $p$ using $0^\jump$ as an oracle.  Consider the algorithm that computes approximations to
$0^\jump$ by simultaneously simulating all programs on input 0 and keeping track of which programs have
halted. As these approximations to $0^\jump$ are gradually revealed, the algorithm may use program $p$ with
them in an attempt to produce the real number $z$.  All of the proper approximations to $0^\jump$ are
actually writable. If one of the proper approximations can successfully produce $z$, then $z$ also is
writable, and consequently $0 \iequiv  z$. Conversely, if none of the proper approximations can produce
$z$, then on input $z$ we can recognize $0^\jump$ as the true approximation, the first approximation able
to produce $z$, and consequently $z \iequiv 0^\jump$.  So in no case can the real $z$ be strictly in the
middle.

Nevertheless, despite this negative solution of Post's problem in the context of real degrees, there is a
supertask analogue of the positive answer in the classical theory when it comes to sets of reals.
Specifically, there is an oracle $A$, a set of reals, that is strictly in between 0 and the jump $0^\jump$.
Indeed, there are incomparable semi-decidable sets $A \perp_\infty B$ (meaning that neither set is
computable from an oracle for the other), and a rich structure of degrees between 0 and $0^\jump$.  The
argument generalizes the Friedburg-Munchnik priority argument (see Soare 1987) to the supertask context,
much as Sacks did for $\alpha$-recursion theory (see Sacks 1990).  So the vibrant classical study of
degrees between 0 and the jump $0^\jump$ that I mentioned earlier finds its analogue here, in the context
of the degrees of sets of reals.

Let me close this article by asking the broad questions: What is the structure of the infinite time Turing
degrees? To what extent do the properties of this structure mirror or differ from the classical structure?
The field is wide open.

\QuietSection References

\begin{list}{}{}

\item J. R. Buchi, Logic, Methodology and Philosophy and Science (Proc. 1960 Internat. Congr.), volume 1 of 1,
Stanford University Press, Stanford, CA, 1962.

\item C. S. Chihara, On the possibility of completing an infinite process, Philosophical Review, 74:74-87,
1965.

\item Jack Copeland, Super Turing-Machines, Complexity, vol. 4 (1998a) pp. 30-32.

\item Jack Copeland, Even Turing Machines Can Compute Uncomputable Functions, In Calude, C., Casti, J.,
Dinneen, M. (eds) Unconventional Models of Computation, London: Springer-Verlag, 1998b, pp. 150-164.

\item Jack Copeland, Accelerating Turing Machines, Minds and Machines (in press, to appear in 2002, in a
Special Issue on effective procedures, edited by Carol Cleland).

\item John Earman, Bangs, crunches, whimpers and shrieks: singularities and acausalities in relativistic
spacetimes, The Clarendon Press, Oxford University Press, New York, 1995.

\item John Earman and John D. Norton, Forever is a day: supertasks in Pitowski and Malament-Hogarth
spacetimes, Philos. Sci., 60(1):22-42, 1993.

\item Solomon Feferman and C. Spector, Incompleteness Along Paths in Progressions of Theories, The Journal of
Symbolic Logic, 27:383-390, 1962.

\item Joel David Hamkins and Andy Lewis, Infinite time Turing machines, The Journal of Symbolic Logic,
65(2):567-604, 2000.

\item Joel David Hamkins and Daniel Seabold, Infinite time Turing machines with only one tape, Mathematical
Logic Quarterly, 47(2):271--287, 2001.

\item Joel David Hamkins and Andy Lewis, Post's problem for supertasks has both positive and negative
solutions, to appear in the Archive for Mathematical Logic.

\item Hogarth, Does general relativity allow an observer to view an eternity in a finite time? Foundations of
Physics Letters, 5:173-181, 1992.

\item Hogarth, Non-Turing computers and non-Turing computability, volume 1 of East Lansing: Philosophy of
Science Association (D.  Hull, M.  Forbes and R.  B.  Burian, eds.), 126-138,  1994.

\item Benedikt L\"owe, Revision sequences and computers with an infinite amount of time, Logic Comput.,
11(1):25-40, 2001.

\item J. P. Laraudogoitia, A beautiful supertask, Mind, 105(417):81-84, January 1996.

\item Pitowsky, The physical Church thesis and physical computational complexity, Iyyun, 39:81-99, 1990.

\item Gerald E.  Sacks, Higher Recursion Theory, Springer-Verlag Publishing Company, Berlin, 1990.

\item Robert I. Soare, Recursively Enumerable Sets and Degrees, Springer-Verlag Publishing Company, New York,
1987.

\item Thomson, Tasks and super-tasks, Analysis, XV:1-13, 1954-55.

\item Philip Welch, Friedman's trick: Minimality arguments in the infinite time Turing degrees, in ``Sets and
Proofs", Proceedings ASL Logic Colloquium, 258:425-436, 1999.

\item Philip Welch, The lengths of infinite time Turing machine computations, Bulletin of the London
Mathematical Society, 32(2):129-136, 2000a.

\item Philip Welch, Eventually infinite time Turing machine degrees: Infinite time decidable reals, Journal
of Symbolic Logic, 65, No.3, 2000, pp11, 2000b.

\end{list}

\QuietSection Notes

\small

\begin{enumerate}
\item My research has been supported in part by grants from the CUNY Research Foundation and the National
Science Foundation.

\item I heard a version of this example in an email exchange with Peter Vallentyne in March, 2001.

\item The editor Jack Copeland has pointed out that Church and Turing put forth only the weaker thesis
asserting that the functions computable by effective procedures that could in principle be carried out by a
person with paper and pencil are Turing machine computable.

\item While the binary expansion of a real number reveals a close connection between ordinary real numbers in
the unit interval and infinite binary sequences, they are not entirely equivalent, largely because the
binary expansion of a real number may not be unique.  For example, while $0.0111\cdots$ and $0.1000\cdots$
are equal as real numbers, the sequences $0111\cdots$ and $1000\cdots$ are distinct.  We may regard the set
of natural numbers as naturally included in the space of infinite binary sequences by identifying the
number 0 with the sequence $000\cdots$, the number 1 with $100\cdots$, the number 2 with $110\cdots$, and
so on.

\item The exact complexity class of the decidable sets lies between $\Sigma^1_1 \Union \Pi^1_1$ and
$\Delta^1_2$, since every semi-decidable set has complexity $\Delta^1_2$.  The Arithmetic sets are exactly
those that are uniformly decidable with a finite number of limits, that is, before $\omega^2$, and the
hyperarithmetic sets, the $\Delta^1_1$ sets, are exactly the sets which can be decided in some bounded
recursive ordinal length of time, uniformly before $\omega_1^{CK}$.  The class $\Delta^1_2$ is known to be
closed under the supertask jump operators, so there is a very rich supertask complexity hierarchy entirely
contained within $\Delta^1_2$.

\item Friends with children have informed me that such an expectation is unwarranted; one sometimes can't get
the child to stop at the right time.  This reminds me of a time when my younger brother was in
kindergarten, the children all sat in a big circle taking turns saying the next letter of the alphabet: A,
B, C, and so on, around the circle in the manner of the usual children's song (evidently an American
phenomenon).  After the letter K, the next child contributed LMNOP, thinking that this was only one letter.

\end{enumerate}

\end{document}